\documentclass[preprint,review,superscriptaddress]{amsart}
\usepackage{times}
\usepackage{amssymb}
\usepackage{dcolumn}
\usepackage{bm}
\usepackage{amsfonts,latexsym,graphicx,verbatim}
\usepackage[matrix,frame,arrow]{xy}
\usepackage{amsmath}
\usepackage{subfigure}
\newcommand{\OO}{\mathrm{O}}
\newcommand{\dsp}{\displaystyle}

\usepackage{amsthm}

\begin{document}
\title{Evolution of Weak Shocks in One Dimensional Planar and Non-planar Gasdynamic Flows}
\author{V D Sharma}
\address{Department of Mathematics, Indian Institute of Technology Bombay.}
\author{Raghavendra Venkatraman}
\address{Department of Mechanical Engineering, Indian Institute of Technology Roorkee.}
\email{raghav.16venkat@gmail.com}
\begin{abstract}
Asymptotic decay laws for planar and nonplanar shock waves and the first
order associated discontinuities that catch up with the shock from behind
are obtained using four different approximation methods. The singular
surface theory is used to derive a pair of transport equations for the
shock strength and the associated first order discontinuity, which
represents the effect of precursor disturbances that overtake the shock
from behind. The asymptotic behaviour of both the discontinuities is
completely analysed. It is noticed that the decay of a first order
discontinuity is much faster than the decay of the shock; indeed, if the
amplitude of the accompanying discontinuity is small then the shock decays
faster as compared to the case when the amplitude of the first order
discontinuity is finite (not necessarily small). It is shown that for a
weak shock, the precursor disturbance evolves like an acceleration wave at
the leading order. We show that the asymptotic decay laws for weak shocks
and the accompanying first order discontinuity are exactly the ones
obtained by using the theory of nonlinear geometrical optics, the theory
of simple waves using Riemann invariants, and the theory of relatively
undistorted waves. It follows that the relatively undistorted wave
approximation is a consequence of the simple wave formalism using Riemann
invariants.
\end{abstract}

\maketitle
\section{Introduction}
The study of the evolutionary behaviour of nonlinear waves in diverse branches of continuum mechanics has long been a subject of great interest from both mathematical and physical points of view. A rigorous mathematical approach to describe the kinematics of a shock has been proposed by Maslov \cite{maslov} using the theory of generalized functions. General properties and the instantaneous growth and decay behaviour of shock waves in elastic materials have been studied using the singular surface theory in a series of articles by Chen \cite{chen}, Nunziato and Herrmann \cite{nunziato} , Chen and Gurtin \cite{chengurtin}, Eringen and Suhubi \cite{suhubi}, McCarthy \cite{McCarthy}, Ting \cite{ting}, Bailey and Chen \cite{bailey} and, Quintanilla and Straughan \cite{quintanilla}. Evolutionary behaviour of shocks in fluids, using singular surface theory, has been discussed by Grinfeld \cite{grinfeld}, Anile \cite{anileprop}, Straughan \cite{straughan1, straughan2}, Jordan \cite{jordan1,jordan2}, Radha et al\cite{vds02}, and Pandey et al \cite{vds09}. Using a procedure based on the kinematics of one-dimensional motion, Sharma and Radha \cite{vds1994,vds3D} and Batt and Ravindran \cite{Batt} have studied the evolutionary behaviour of shocks. Numerical approaches such as the Godunov method have also been used to understand shock evolutionary behaviour, see for instance Lax \cite{Lax}, Randall LeVeque \cite{leveque1} and Christov and Jordan \cite{na}. For a detailed exposition on the Godunov method and its extensions, the reader is referred to \cite{leveque2}. An interesting study on the evolutionary behaviour of shocks in elastic nonconductors has been carried out by Fu and Scott \cite{Fu2,Fu3} using singular surface theory together with the method of multiple scales, method of weakly nonlinear geometrical optics and the shock fitting method.  

The essential approximation involved in the method of weakly nonlinear geometrical optics is that the wave is of small amplitude and is slowly modulated; this means that the wave profile changes shape over distances or times much larger than a typical distance or time over which the wave profile itself varies. This method, unlike the method of characteristics, is applicable for problems involving several independent variables; indeed, a formal and systematic derivation of nonplanar simple wave solutions is possible only in the high frequency or geometrical-acoustics limit. Modulated simple wave theories have been developed by several authors; papers by Parker \cite{parker1,parker2} and an extensive review by Seymour and Mortell \cite{seymour} are particularly illuminating. For the generalization of the weakly nonlinear theory in high frequency approximation and for applications of this work, we refer the reader to the work carried out by Hunter and Keller \cite{Hunter}, Hunter \cite{Hunter1}, Anile et al \cite{Anile3}, Srinivasan and Sharma \cite{GKSVDS} and Ali and Hunter \cite{ali}. 

The other approximation method, namely the simple wave formulation using Riemann invariants is useful in the construction of solutions and development of singularities. For hyperbolic systems involving source terms or non-planar geometry, the governing equations may not admit exact simple wave solutions, but under small wave amplitude assumption, they may admit asymptotic modulated simple wave solutions. Indeed, hyperbolic systems are mathematically tractable using simple waves and Riemann techniques of gasdynamics, see for instance, Zheng \cite{zheng}, Sterck \cite{sterck} and Li et al \cite{li}. Lastly, we consider the method of relatively undistorted waves; this applies to solutions of hyperbolic partial differential equations, and is based on the same slow modulation approximation used in the method of weakly nonlinear geometrical optics. But, unlike the method of weakly nonlinear geometrical optics, this theory makes no assumption on the magnitude of a disturbance (see, Varley and Cumberbatch \cite{varley,varley1} and Seymour and Mortell \cite{seymour}).

It may be recalled that an analytic description of the
evolutionary history of a shock of arbitrary strength, both
in fluids and solids, still remains an open problem. However,
kinematics of a weak shock has been studied in the references
cited as above, mostly in solids, using singular surface
theory and shock fitting methods. The purpose of the present
paper is to study planar, axially and radially symmetric
flows of a polytropic gas behind a weak shock, propagating
into a uniform state at rest, using four different
approximation methods and establish their equivalence
within their common range of validity. To the best of
our knowledge, this forms not only the first comprehensive
study of one-dimensional planar and non-planar weak shocks
in ideal gasdynamic flows, but also provides a comparative
analysis of four seemingly different approaches. The approximation methods employed herein for the derivation
of asymptotic decay laws for weak shocks and the
accompanying first order discontinuities are
i) the singular surface theory, ii) the method of
weakly nonlinear geometrical optics, iii) the characteristic
approach using simple waves as kinematic waves, and
iv) the method of relatively undistorted waves. 

The paper is organized as follows: after setting up the basic
equations and jump conditions across a shock of arbitrary
strength, we devote sections 3 and 4 to the derivation of
transport equations for the variation of
jumps in pressure and its space derivatives; the equations
are coupled with those involving jumps in higher order space
derivatives of pressure. In section 5, we analyze the pair
of coupled transport equations for the shock strength and
accompanying first and second order discontinuities; the shock
strength is assumed to be small but the associated first order
jump discontinuity may be finite or small. Under these
assumptions, the transport equations are solved exactly to
the leading order once the strength of the accompanying
second order discontinuity is neglected, and the evolutionary
behaviour of both the shock and the associated first order
discontinuity is completely analysed. Furthermore, the transport equations have also been numerically integrated in order to elucidate the results of this section. In section 6, we show how
the method of weakly nonlinear geometrical optics leads to the
same evolution laws for the shocks and first order
discontinuities as those obtained in the preceding section.
Sections 7 and 8 deal with the other two approximation methods,
where it is shown that they lead to the same evolution laws
within their common range of validity. The main results and
conclusions inferred from this work are given in the last section.

\section{Basic Equations and Jump Conditions across Shocks}
We consider the one dimensional unsteady flow of a polytropic gas behind a planar, cylindrical or a spherical  wave propagating into a uniform state at rest.The equations governing the flow are \cite{whitham}
\begin{equation}
\label{1}
\begin{array}{l}
\rho_t + u \rho_x + \rho u_x + \frac{j}{x}\rho u = 0; \hspace{0.2cm} u_t + uu_x + \frac{1}{\rho}p_x = 0; \hspace{0.2cm} p_t + up_x + \rho a^2\left(u_x + \frac{ju}{x}\right) = 0, 
\end{array}
\end{equation}
where $\rho$ is the density, $u$ the gas velocity, $p$ the pressure, and $a = (\gamma p/\rho)^{1/2}$ the sound speed with $\gamma$ as the specific heat ratio of the gas; here $t$ stands for time and $x$ the distance being either axial in flows with planar $(j=0)$ geometry, or radial in cylindrically symmetric $(j=1)$ and spherically symmetric $(j=2)$ flow configurations.

As system \eqref{1} is quasilinear with each of these equations being a direct consequence of the corresponding conservation law, one expects shocks to appear in the flow after a finite running length or time. We consider the gas motion containing a shock wave propagating into a homogeneous quiescent equilibrium gas with intrinsic velocity $U$ given by $U^{-1} = s'(x) = \tfrac{\,dt}{\,dx}$ where $t = s(x)$ denotes the location of the shock at position $x.$ 
Across the shock front, we have the following Rankine Hugoniot (R-H) conditions,
\begin{equation} 
\label{2}
\begin{array}{l}
U[\rho] = [\rho u]; \hspace{0.4cm} U[\rho u] = [p + \rho u^2]; \hspace{0.4cm}U[ \frac{1}{2} \rho u^2 + \rho e] = [\left(\frac{1}{2} \rho u^2 + \rho e\right)u + pu], 
\end{array}
\end{equation}
where $[f] = f_- - f_+$ denotes the jump in a variable $f$ with $f_+$ and $f_-$ being the values of $f$ immediately ahead of and immediately behind the shock, and  $e$ is the specific internal energy given by $e = \tfrac{1}{\gamma -1}\tfrac{p}{\rho}.$

In addition, following the singular surface theory, we have the following compatibility condition \cite{truesdell}, 
\begin{equation} \label{3}
\frac{\mathrm{d}[f]}{\mathrm{d}t} = [f_t] + U[f_x],
\end{equation}
where $\tfrac{\,d}{\,dt}$ denotes the time derivative following the shock front. 
The medium ahead of the shock front is assumed to be uniform and at rest, i.e., $\rho_+, p_+$ and $a_+$ are constant, and $u_+ = 0.$ The velocities $u,a$ and $U$ appearing in \eqref{1} and \eqref{2} are nondimensionalized by $a_+,$ and the remaining variables $\rho, p,x$ and $t$ are rendered dimensionless by $\rho_+, \rho_+ a_+^2, x_0, $ and $\tfrac{x_0}{a_+},$ respectively; here $x_0$ characterizes the reference length of the medium. Thus, the variables appearing in \eqref{1} and \eqref{2} will henceforth be regarded as dimensionless; indeed, in the medium ahead of the shock, we have $a = a_+ = 1$ and $\rho = \rho_+ = 1.$ The following expressions follow readily from conditions \eqref{2}, i.e.,
\begin{equation} \label{RH}
\begin{array}{l} 
\dsp [u] = \frac{U[\rho]}{1 + [\rho]}, \hspace{1cm} [p] = U[u],\hspace{1cm} [u] = \frac{2}{\gamma+1}\frac{U^2-1}{U}.
\end{array}
\end{equation}
In view of \eqref{RH}, it follows that, 
\begin{equation} \label{derivatives}
\frac{\,d[p]}{\,dt} = \frac{2U^3}{(U^2+1)}\left(\frac{\,d[u]}{\,dt} \right) = \left(\frac{\mu}{(\gamma+1)}\right)^2 \cdot \left( \frac{\,d[\rho]}{\,dt}\right),
\end{equation}
where $\mu = 2 + (\gamma -1)U^2.$
\section{Transport Equation for Shock Strength using Singular Surface Theory}
In order to derive the equation for the shock strength, we take jumps in the Euler equations \eqref{1}, and using the conditions \eqref{3}-\eqref{derivatives} to obtain 
\begin{equation} \label{jumptransport}
\frac{\,d[\mathbf{V}]}{\,dt} + \mathcal{A}[\mathbf{V_x}] + \Omega [u]\mathcal{B}= 0,
\end{equation}
where, 
\begin{subequations} 
\begin{eqnarray} \nonumber
[\mathbf{V}] = \left( \begin{array}{ccc}
\dsp [\rho] & [u] &  [p]
\end{array}\right)^{tr}, \hspace{1cm} [\mathbf{V_x}] = \left( \begin{array}{ccc}
\dsp [\rho_x] &
 [u_x] &
 [p_x]
\end{array}\right)^{tr}, \\  \nonumber
\mathcal{A} = \left( \begin{array}{ccc}
[u]-U & 1 + [\rho] & 0\\
0 & [u]-U & \frac{1}{1+[\rho]}\\
0 & 1 + \gamma [p] & [u] - U
\end{array}\right), \hspace{1cm}
\mathcal{B} = \left(\begin{array}{l}
1 + [\rho]\\
0 \\
1 + \gamma [p]
\end{array}\right),
\end{eqnarray}
\end{subequations}
and $\Omega = \tfrac{j}{x(t)}$ . It may be noticed that the form \eqref{jumptransport} is not only convenient for algebraic manipulation, but is also helpful in understanding the pattern of transport equations for higher order discontinuities.

Notice that \eqref{jumptransport} contains the unknown quantities $[u_x], [p_x]$ and $[\rho_x];$ eliminating between these equations the unknown discontinuities $[\rho_x]$ and $[u_x]$ by taking a suitable linear combination, we arrive at the first transport equation governing the shock strength
\begin{equation} \label{main1}
\frac{\,d[p]}{\,dx} = k_{11}[p_x] + k_{12},
\end{equation}
where the coefficients $k_{11}$and $k_{12}$ are given by
\begin{subequations} \label{coeff1}
\begin{eqnarray} \nonumber
k_{11} = -\frac{2(U^2-1)\mu}{U^2(2\mu+\nu)+ \nu},\hspace{0.6cm} k_{12} = -\frac{4(U^2-1)\mu\nu}{U^2(2\mu+\nu)+ \nu}\left( \frac{\Omega}{(\gamma+1)^2}\right),
\end{eqnarray}
\end{subequations}
with $\nu = 2\gamma U^2 + 1 - \gamma.$
We note that an immediate consequence of \eqref{main1} cannot reveal the complete history of the evolutionary behaviour of the shocks under consideration because of the appearance of an unknown entity $[p_x]$. It is clear from \eqref{main1} that the evolutionary behaviour of shocks at any time $t$ depends not only on the strength of the shock, its curvature $\tfrac{j}{x(t)}$ and the specific heat ratio $\gamma,$ but also on the pressure gradient  immediately behind the wave $[p_x]$. In the following section we determine a transport equation for the unknown $[p_x],$ which represents the effect of disturbances that overtake the shock from behind. 

\section{Transport Equation for the Unknown term $[p_x]$}
As the coupling term $[p_x]$ in \eqref{main1} is unknown we need to obtain a transport equation for it; to this end, we first note that the $[u_x]$ and $[p_x]$ eliminant of equations \eqref{jumptransport} along with \eqref{main1} yields a relation between $[\rho_x]$ and $[p_x]$. Similarly the $[\rho_x]$ and $[p_x]$ eliminant of equations \eqref{jumptransport} along with \eqref{main1} yields a relation between $[u_x]$ and $[p_x]$. In this section, we determine the transport equation for the unknown pressure gradient jump as described in the preceding section. We first note that from equations \eqref{jumptransport}
\begin{equation} \label{relations}
\left( \begin{array}{l}
\dsp [u_x]\\
\dsp [\rho_x]
\end{array}\right) = \mathcal{T}\left( \begin{array}{l}
\dsp [p_x]\\
\dsp 1
\end{array}\right),
\end{equation}
where $\mathcal{T} = [T_{ij}]_{2 \times 2}$ is defined by
\begin{equation}
\mathcal{T} = \left(\begin{array}{cc}
\frac{\mu - (\gamma + 1)k_{11} U^2}{\nu U} & -\frac{k_{12} (\gamma+1)(U^4-1)}{k_{11}U((2\mu+\nu)U^2 + \nu)} \\
\frac{(\gamma+1)^2 U (\mu^2 U - (\gamma+1)(\mu U^2 - \nu)k_{11})}{\nu \mu^3} & \frac{U(\gamma+1) ((\gamma+1)^2 (U^2 + 3)k_{12} + 4 \mu U(U^2-1)\Omega)}{2\mu^3}
\end{array}\right).
\end{equation} 

Under this setting, we now differentiate equations \eqref{1}, take jumps across the shock, and use the shock conditions \eqref{3} and \eqref{RH} to obtain the following system of equations
\begin{equation} \label{jumpdisc}
\frac{\,d[\mathbf{V_x}]}{\,dt} + \mathcal{A}[\mathbf{V_{xx}]} + \Omega \mathcal{D}[\mathbf{V_x}] + \Omega' [u]\mathcal{B} + \mathcal{C} = 0,
\end{equation}
where $[\mathbf{V}], [\mathbf{V_x}], \mathcal{A},$ and $\mathcal{B} $ have the same meaning as before and 
\begin{equation}
\nonumber
\begin{array}{l}
[\mathbf{V_{xx}}] = \big( \begin{array}{ccc} [\rho_{xx}] & [u_{xx}] & [p_{xx}]\end{array}\big)^{tr}, \\
\mathcal{D} = \left(\begin{array}{ccc}
[u] & (1+ [\rho]) & 0 \\
0 & 0 & 0\\
0 & (1 + \gamma [p]) & \gamma(\gamma-1)[u]
\end{array} \right), \hspace{1cm} \mathcal{C} = \left(\begin{array}{l}
2[u_x][\rho_x] \\
- \left( [u_x]^2 + \frac{[\rho_x][p_x]}{(1+[\rho])^2}\right) \\
(\gamma+1)[u_x][p_x]
\end{array}\right).
\end{array}
\end{equation}
Eliminating the unknown discontinuity terms $[\rho_{xx}]$ and $[u_{xx}]$ between equations \eqref{jumpdisc}, and using \eqref{2} and \eqref{relations} we arrive at the second transport equation that governs the discontinuity $[p_x]$, i.e.,
\begin{equation} \label{main2}
\frac{\,d[p_x]}{\,dx} + k_{21} [p_{xx}] + k_{22} [p_x]^2 + k_{23}[p_x] + k_{24} = 0,
\end{equation}
where 
\begin{equation}
\begin{array}{l}
\dsp k_{21} = \frac{(U^2-1)\eta}{U^2}, \hspace{1cm} \eta = \frac{\mu}{2\mu-(\gamma+1)Uk_{11}}, \\\\
\dsp k_{22} = \left( \frac{(\gamma+1)\eta}{U\mu} \right) \left( T_{11} \left(\mu + \frac{\nu U T_{11}}{(\gamma+1)}\right) + \left( \frac{\nu k_{11}}{4}\right) \left(\frac{\,dT_{11}}{\,dU}\right)\right) - \frac{\mu \nu\eta T_{21}}{(\gamma+1)^2U^4}, \\\\
\dsp k_{23} = \left( \frac{T_{12} \eta}{\mu} \right) \cdot \frac{\mu(\gamma+1) + 2\nu U T_{11}}{U} + \eta \Omega \frac{\gamma+1}{U}\left( \nu T_{11}  + \left( \frac{2\gamma}{U}\right)(U^2-1)\right) \\ \dsp \phantom{k_{23} =}- \frac{\mu\nu\eta T_{22}}{U^4(\gamma+1)^2} + \left(\frac{\eta \nu k_{11}}{4\mu}\right)\frac{\gamma+1}{U}\left( \frac{k_{12}}{k_{11}}\frac{\,dT_{11}}{\,dU} + \partial_U T_{12}|_x\right), \\\\
\dsp k_{24} = 2\eta \left( \frac{\nu}{U^2}\right)(U^2-1)\frac{\Omega'}{(\gamma+1)^2} + \frac{\eta \nu T_{12} \Omega}{U(\gamma+1)} \\ \phantom{k_{24} =} + \left( \frac{\nu U \eta}{U\mu}\right)\left( T_{12}^2 + U(\partial_x T_{12})\big|_U + (\gamma+1)\left(\frac{k_{12}}{4U}\right) \partial_U T_{12}\big|_x\right). \\
\end{array}
\end{equation}
It may be noticed from \eqref{main2} that an analytical description of $[p_x]$ is again obscured by the presence of higher order jump discontinuity $[p_{xx}]$ which is an unknown, whose determination requires a repetition of the above procedure; proceeding in this manner we arrive at an open system consisting of an infinite set of transport equations for the coupling terms - the jumps in higher order derivatives of $p$. In order to provide a natural closure on the hierarchy of these equations, we consider the coupled system \eqref{main1} and \eqref{main2} and render it tractable by making assumptions on the higher order discontinuity $[p_{xx}]$ (see Fu and Scott\cite{Fu2}). 
\section{Evolution Laws for Weak Shocks}
For a weak shock we may assume that $[p] = \OO(\epsilon), 0 < \epsilon << 1,$ and consider the two cases, $[p_x] = \OO(1)$ and $[p_x] = \OO(\epsilon).$ In the former case, direct perturbation ansatz yields a uniformly valid solution to the leading order approximation, while in the second case we use the method of multiple scales. We now present these details below.
 
For a weak shock, it follows from \eqref{RH} that 
\begin{equation} \label{15star}
\dsp [u] = [p] + \OO(\epsilon^2); \hspace{0.6cm} [\rho] = [p] + \OO(\epsilon^2); \hspace{0.6cm} U = \left(1 + \frac{\gamma +1}{4}[p]\right) + \OO(\epsilon^2),
\end{equation} 
and thus the equations \eqref{main1} and \eqref{main2} can be approximated to give 
\begin{equation} \label{15}
\begin{array}{l}
\frac{\,d[p]}{\,dx} + \frac{\gamma+1}{4}[p][p_x] + \frac{j}{2x}[p] = 0, \\
\frac{\,d[p_x]}{\,dx} + \frac{\gamma+1}{4}[p][p_{xx}] + \frac{\gamma+1}{2}[p_x]^2 + \frac{j}{2x}[p_x] = 0.
\end{array}
\end{equation}
We need to solve these equations assuming that the initial conditions are $[p]|_{x=1} = h$ and $[p_x]_{x=1} = k,$ where $|h|$ is assumed to be small and can be taken as $\epsilon.$

\subsection{Case 1: $[p_x] = \OO(1)$}
When the effect of the first order induced discontinuity or the disturbances that overtake the shock from behind is strong, i.e., $[p_x] = \OO(1),$ equations \eqref{15} together with $[p] = \epsilon Y(x), [p_x] = Z(x)$ and $[p_{xx}] = \OO(\epsilon^q), q \geq 0$ yield to the leading order the following system
\begin{equation}\label{system1}
\frac{\,dY}{\,dx} + \frac{\gamma+1}{4}YZ + \left(\frac{j}{2x}\right) Y = 0, \hspace{0.6cm} \frac{\,dZ}{\,dx} + \frac{\gamma+1}{2}Z^2 + \left( \frac{j}{2x}\right) Z =0,
\end{equation}
with the initial conditions 
\begin{equation} \label{IC1}
Y(1) = \tfrac{h}{\epsilon}, \hspace{1cm} Z(1) = k.
\end{equation}
Equation $(16b)$ is the well known Bernoulli equation that governs the amplitude $[p_x]$ of an acceleration wave; a complete analysis of such an equation modifying and generalizing several known results in the literature has been reported in \cite{JMP1981} and \cite{AA1983}. We thus conclude that the precursor disturbance that overtakes the shock from behind evolves like an acceleration wave at the leading order. System \eqref{system1}-\eqref{IC1} can be directly integrated to yield
\begin{equation}\label{integral1}
Y = h I^{-\frac{1}{2}}\psi(x), \hspace{1cm} Z = k I^{-1} \psi(x),
\end{equation}
where $\psi(x) = x^{-\frac{j}{2}}$ and
\begin{equation}
I = 1 + \frac{\gamma+1}{2}k J(x), \hspace{1cm} J(x) = \int_1^x \psi(s) \,ds.
\end{equation}
It may be noticed that if $k > 0,$ equations \eqref{integral1} imply the following asymptotic decay laws for the shock and the associated first order discontinuity  
\begin{equation} \label{p1anal}
\dsp [p] \sim h \left( \frac{2}{(\gamma +1)k}\right)^{\tfrac{1}{2}} \left\{ 
\begin{array}{cc}
x^{-1/2}, & \mbox{ plane }, \\
\tfrac{1}{\sqrt{2}}x^{-3/4}, & \mbox{ cylindrical},\\
x^{-1}(\log x)^{-1/2}, & \mbox{spherical};
\end{array}\right.
\end{equation}
and 
\begin{equation} \label{px1anal}
\dsp [p_x] \sim h \left( \frac{2}{\gamma +1}\right) \left\{ 
\begin{array}{cc}
x^{-1}, & \mbox{ plane }, \\
\tfrac{1}{2}x^{-1}, & \mbox{ cylindrical},\\
(x \log x)^{-1}, & \mbox{spherical},
\end{array} \right.
\end{equation}
as $ x \to \infty;$ the asymptotic law for the shock is in full agreement with earlier results obtained using simple wave theory \cite{whitham,vdsbook}. It may be noticed that the decay of the induced discontinuity $[p_x]$ is much faster than the decay of the shock. 
 \begin{figure}[ht] \label{fig1} \centering
\includegraphics[scale=.5]{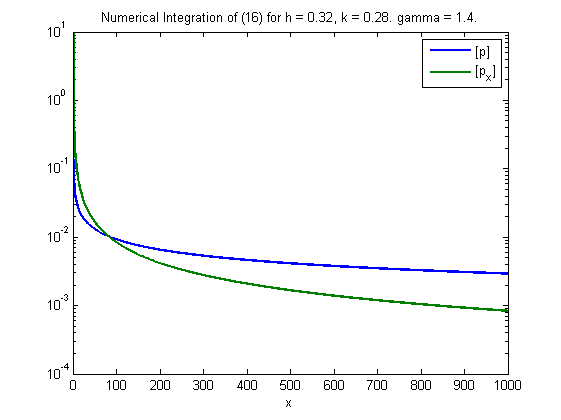}
\caption{Numerical solution of system \eqref{15}: Planar Configuration  $h = 0.32, k = 10, \gamma = 1.4$}
\end{figure}
In order to further elucidate the results of this section, the transport equations \eqref{15} were numerically integrated for the planar case for an ideal gas with $\gamma = 1.4$. Figure 1 shows the results of the numerical integration. Table 1 further establishes the validity of asymptotic results \eqref{p1anal} and \eqref{px1anal} by comparing them with results of numerical integration of system \eqref{15}.

\begin{table}
\centering
\label{table1}
\scalebox{0.8}{\begin{tabular}{|c|c|c|c|c|}
\hline
&\multicolumn{2}{c}{$h= 0.32, k = 10$}&
\multicolumn{2}{c}{$h=0.32, k = 0.28$}\\
\cline{1-5}
Position $x$ & Error in $[p]$ & Error in $[p_x]$&Error in $[p]$ & Error in $[p_x]$\\
\hline \hline
1.476 &	4.332e-2	& 9.545e-1 &3.374e-2 & 6.752e-2 \\
4.565  &	4.827e-3 & 4.914e-2& 1.317e-2&1.885e-2\\
7.668 &	2.076e-3	&1.593e-2& 7.448e-3&8.723e-3\\
9.563 &	1.464e-3 &9.990e-3& 5.675e-3& 6.133e-3\\
13.3 &	8.752e-4 &	5.032e-3&3.711e-3 & 3.482e-3\\
27.95 &	2.798e-4 &	1.100e-3&1.373e-3& 9.242e-4\\
45.57	&1.332e-4	&4.088e-4& 6.879e-4&3.678e-4\\
65.31	&7.732e-5	&1.979e-4& 4.078e-4&1.832e-4\\
76.04   &6.145e-5      &1.457e-4& 3.271e-4&1.366e-4\\
86.34	&5.075e-5	&1.129e-4&2.716e-4&1.065e-4\\
96.35	&4.301e-5	&9.054e-5&2.311e-4&8.591e-5\\
99.95	&4.07e-5		&	8.411e-5&2.205e-4 & 8.072e-5\\
100	& 4.067e-5		&	8.403e-5&2.190e-4&7.994e-5\\
\hline
\end{tabular}
}
\caption{Numerical and asymptotic results for shock strength and accompanying discontinuity}
\end{table}

\subsection{Case 2: $[p_x] = \OO(\epsilon)$} 
When the induced discontinuity $[p_x]$ is weak,i.e., $[p_x] = \OO(\epsilon),$ \eqref{15} together with $[p] = \epsilon Y(x),[p_x] = \epsilon Z(x)$ and $[p_{xx}] = \OO(\epsilon^q),$ where $q >1,$ yields to next order the following system
\begin{equation} \label{22}
\frac{\,dY}{\,dx} + \frac{\gamma +1}{4}\epsilon YZ + \left( \frac{j}{2x}\right)Y = \OO(\epsilon^2), \hspace{0.2cm}
\frac{\,dZ}{\,dx} + \frac{\gamma +1}{2}\epsilon Z^2 + \left( \frac{j}{2x}\right)Z = \OO(\epsilon^2).
\end{equation}
We notice that the straightforward perturbation ansatzes of the forms $Y(x) = Y_0(x) + \epsilon Y_1 (x) + \OO(\epsilon^2) ,$  and $Z(x) = Z_0(x) + \epsilon Z_1(x)+ \OO(\epsilon^2) $ yield secular terms that blow up as $x \to \infty$ and render the solutions nonuniformly valid. In order to circumvent this situation, we use the method of multiple scales, which requires the introduction of slow variables \cite{nayfeh}. Now we take up the following cases separately. 
\subsubsection*{Planar $(j=0)$ waves:} 
In this case, we introduce the slow variable, $\eta = \epsilon x,$ and introduce the formal expansions
\begin{equation}
\nonumber
Y = Y_0(x,\eta) + \epsilon Y_1(x,\eta) + \OO(\epsilon^2),\hspace{0.6cm}
Z = Z_0(x,\eta) + \epsilon Z_1(x,\eta) + \OO(\epsilon^2), 
\end{equation} 
which on using in \eqref{22} yield
\begin{equation} \label{plan1}
\begin{array}{cc}
\OO(1): & \frac{\partial Y_0}{\partial x}  = 0, \hspace{1cm} \frac{\partial Z_0}{\partial x}  =0, \\
\OO(\epsilon) : & \frac{\partial Y_1}{\partial x} + \frac{\partial Y_0}{\partial \eta} + \frac{\gamma + 1}{4}Y_0Z_0  = 0, \hspace{0.6cm} \frac{\partial Z_1}{\partial x} + \frac{\partial Z_0}{\partial \eta} + \frac{\gamma+1}{2}Z_0^2  =0,
\end{array}
\end{equation}
with the following initial conditions
\begin{equation} \label{bc1}
Y_0(1) = \tfrac{h}{\epsilon}, \hspace{1cm} Z_0(1) = \tfrac{k}{\epsilon},
\end{equation}
where $|h|$ and $|k|$ are small.
The first pair of equations in \eqref{plan1} yield
\begin{equation}
Y_0 = f(\eta) , \hspace{1cm} Z_0 = g(\eta) .
\end{equation}
Furthermore, in order to avoid secular terms from the equations at this order $\OO(\epsilon)$, we have 
\begin{equation}
\label{sec1}
f'(\eta) + \frac{\gamma+1}{4}f(\eta)g(\eta) = 0, \hspace{1cm} g'(\eta) + \frac{\gamma+1}{2}g^2(\eta) = 0.
\end{equation} 
Finally solving \eqref{sec1} for $k > 0$ in conjunction with the initial conditions \eqref{bc1}, we obtain
\begin{equation} \label{panal}
[p] \sim h \left(\frac{2}{(\gamma +1)k}\right)^{1/2}x^{-\frac{1}{2}} , \hspace{1cm} [p_x] \sim \frac{2}{(\gamma +1)}x^{-1} \hspace{1cm} \mbox{as } x \to \infty.
\end{equation}
\begin{figure}[ht] \label{fig2} \centering
\includegraphics[scale=.5]{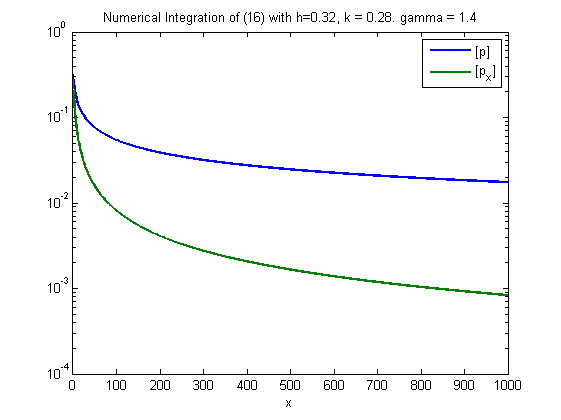}
\caption{Numerical solution of system \eqref{15}: Planar Configuration $h = 0.32, k = 0.28, \gamma = 1.4$}
\end{figure}
In order to further elucidate the results of this section, the transport equations \eqref{15} were numerically integrated for an ideal gas with $\gamma = 1.4$ Figure 2 shows the results of the numerical integration of \eqref{15}. Table 1 further establishes the validity of the asymptotic results in \eqref{panal} with the results of numerical integration of \eqref{15}.
It may also be noticed that in this case, the shock decays faster as compared to the previous case when $k = \OO(1);$ this has been shown in Figure 3. 
\begin{figure}[ht] \label{fig3} \centering
\includegraphics[scale=.5]{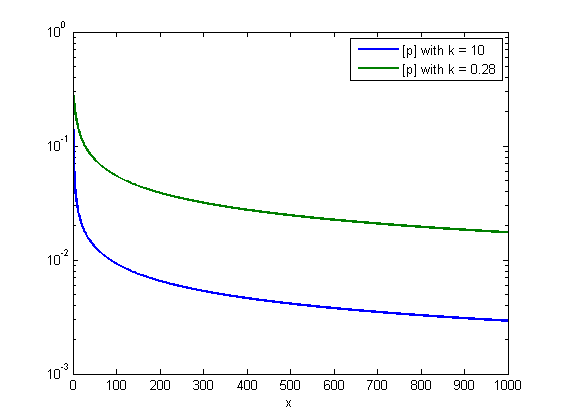}
\caption{Comparison of shock strength decay in case 1 and case 2 above}
\end{figure}
Finally, it is interesting from \eqref{panal} and \eqref{px1anal} that the asymptotic behaviour of the precursor disturbance $[p_x]$ is independent of $k,$ this has been illustrated in Figure 4 where the two curves are indistinguishable in the asymptotic limit. 
\begin{figure}[ht] \label{fig3} \centering
\includegraphics[scale=.5]{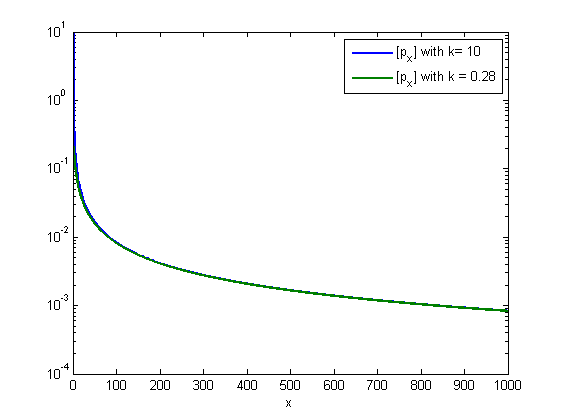}
\caption{Comparison of $[p_x]$ in case 1 and case 2 above}
\end{figure}
\subsubsection*{Cylindrical $(j=1)$ wave:}
In this case, the system of equations \eqref{22} with the slow variable $\eta = \epsilon \sqrt{x}$ yields the following hierarchy of partial differential equations
\begin{equation}
\begin{array}{cc}
\OO(1): & \frac{\partial Y_0}{\partial x} + \left( \frac{1}{2x}\right) Y_0 = 0, \hspace{0.6cm} \frac{\partial Z_0}{\partial x} + \left(  \frac{1}{2x}\right) Z_0 = 0, \\
\OO(\epsilon): & \frac{\partial Y_1}{\partial x} + \frac{1}{2\sqrt{x}}\frac{\partial Y_0}{\partial \eta} + \frac{\gamma+1}{4}Y_0Z_0 + \left(  \frac{1}{2x}\right) Y_1 = 0, \\
& \frac{\partial Z_1}{\partial x} + \frac{1}{2\sqrt{x}}\frac{\partial Z_0}{\partial \eta} + \frac{\gamma+1}{2}Z_0^2 + \left(  \frac{1}{2x}\right) Z_1 = 0, \\
\end{array}
\end{equation}
together with the initial conditions \eqref{bc1}. We find that to the leading order $Y_0 = A(\eta) x^{-1/2}$ and $Z_0 = B(\eta) x^{-1/2},$ where $A$ and $B$ satisfy the differential equations 
\begin{equation}
\nonumber 
A' + \frac{\gamma +1}{2}AB = 0, \hspace{1cm} B' + (\gamma+1) B^2 = 0.
\end{equation}
Here prime denotes differentiation with respect to the variable $\eta;$ solving these equations and using \eqref{bc1}, it follows that for $k>0,$ the shock and induced discontinuity decay asymptotically as 
\begin{equation}
\dsp [p] \sim h\left( \frac{1}{k(\gamma +1)}\right)^{1/2}  x^{-3/4}, \hspace{1cm} [p_x] \sim \frac{1}{(\gamma +1)}x^{-1} \hspace{1cm} \mbox{ as } x \to \infty.
\end{equation}
\subsubsection*{Spherical $(j=2)$ waves:}
We repeat the above procedure with the slow variable $\eta = \epsilon \log x,$  appropriate to this case, and find that for $k >0,$ the decay behaviour to the leading order is
\begin{equation}
\dsp [p] \sim h\left( \frac{2}{k(\gamma +1)}\right)^{1/2}  x^{-1} (\log x)^{-1/2}, \hspace{0.6cm} [p_x] \sim \frac{2}{(\gamma +1)}x^{-1}(\log x)^{-1} \hspace{0.6cm} \mbox{ as } x \to \infty.
\end{equation}

\section{Weakly Nonlinear Geometrical Optics}
In this section, we use the theory of weakly nonlinear geometrical optics, which allows waves to propagate on different families of characteristics, to treat shock formation, and to include interaction of the flow with the shock. This is an asymptotic method whose objective is to understand the laws governing the propagation and interaction of small amplitude high frequency waves \cite{Hunter}. Here we show that this method leads to the same evolution laws for shocks and first order discontinuities as those obtained in the preceding section. 

We assume that the disturbance is of small amplitude as in the previous section; the parameter $\epsilon$ used in the previous section suggests to form an asymptotic analysis. We take $(x,t,\theta)$ as independent variables where $\theta = \tfrac{\phi(x,t)}{\epsilon}$ is the fast variable with $\phi(x,t)$ as the phase function to be determined, and $\phi=$ constant is the characteristic wavelet which left the boundary $x=1$ at $\theta = \tfrac{t}{\epsilon}.$ We now seek small amplitude high frequency solutions to the system \eqref{1} of the form
\begin{subequations} \label{ansatzmod}
\begin{eqnarray}
\rho = 1 + \epsilon \rho^{(1)}(x,t,\theta) + \epsilon^2 \rho^{(2)} (x,t,\theta) + \OO(\epsilon^3) ,\\
u = \epsilon u^{(1)}(x,t,\theta) + \epsilon^2 u^{(2)} (x,t,\theta) + \OO(\epsilon^3),\\
p = \tfrac{1}{\gamma} + \epsilon p^{(1)}(x,t,\theta) + \epsilon^2 p^{(2)} (x,t,\theta) + \OO(\epsilon^3). 
\end{eqnarray}
\end{subequations}
where $\rho^{(1)}, u^{(1)}, p^{(1)}$ etc are regular and bounded functions of their arguments together with their derivatives as $\epsilon \to 0;$ using \eqref{ansatzmod} in \eqref{1} along with
\begin{equation}\nonumber
\frac{\partial f}{\partial x}\big|_{t} = \frac{\partial f}{\partial x}\big|_{t,\theta} + \frac{1}{\epsilon}\phi_{x}\frac{\partial f}{\partial \theta}; \hspace{1cm}\frac{\partial f}{\partial t}\big|_{x} = \frac{\partial f}{\partial t}\big|_{x,\theta} + \frac{1}{\epsilon}\phi_{t}\frac{\partial f}{\partial \theta},
\end{equation}
and equating to zero the coefficients of $\epsilon^i, i =0,1,$ we arrive at the following hierarchy of partial differential equations:
\begin{subequations} \label{hier1}
\begin{eqnarray} \label{order1}
\OO(1): \hspace{3.7cm} \phi_x u^{(1)}_\theta + \phi_t \rho^{(1)}_\theta = 0;\hspace{2.5cm}\\
\nonumber \phantom{\OO(1): \hspace{0.8cm}}\phi_x p^{(1)}_\theta + \phi_t u^{(1)}_\theta = 0;\hspace{2.5cm}\\
\nonumber \phantom{\OO(1): \hspace{0.8cm}}\phi_x u^{(1)}_\theta + \phi_t p^{(1)}_\theta = 0, \hspace{2.5cm}\\
\OO(\epsilon): \hspace{0.6cm}\label{ordereps} 
\rho^{(1)}_t + \phi_t \rho^{(2)}_\theta + \phi_x \big((\rho^{(1)}u^{(1)})_\theta + u^{(2)}_\theta \big) + u^{(1)}_x + \tfrac{ju^{(1)}}{x}=0;\hspace{0.2cm}\\
\nonumber \phantom{\OO(\epsilon): \hspace{0.4cm}} u^{(1)}_t + \phi_t u^{(2)}_\theta + \phi_x\big(u^{(1)}u^{(1)}_\theta + p^{(2)}_\theta - \rho^{(1)}p^{(1)}_\theta \big) + p^{(1)}_x =0; \\
\nonumber \phantom{\OO(\epsilon): \hspace{0.6cm}} p^{(1)}_t + \phi_t p^{(2)}_\theta + \phi_x \big( u^{(1)}p^{(1)}_\theta + \gamma p^{(1)}u^{(1)}_\theta + u^{(2)}_\theta \big) + u^{(1)}_x + \tfrac{ju^{(1)}}{x}=0.
\end{eqnarray}
\end{subequations}
Equations \eqref{hier1} admit nontrivial solutions if $\phi$ satisfies the eikonal equation
\begin{equation}
\nonumber
\tfrac{\phi_t}{\phi_x}\left(\big(\tfrac{\phi_t}{\phi_x}\big)^2 -1\right) = 0,
\end{equation}
which shows that the phase function $\phi(x,t)= $ constant represents a propagating wave front of the leading order system associated with \eqref{1}. Indeed, each of the phase speeds, $-\tfrac{\phi_t}{\phi_x},$ which can take the values $-1,0$ and $1,$ gives the speed of propagation of $\phi$ along normal to the wave front. It follows from \eqref{order1} that along the forward facing wave, $-\tfrac{\phi_t}{\phi_x} =1,$ we have 
\begin{equation} \label{alpha}
\rho^{(1)} = u^{(1)}= p^{(1)}= \alpha(x,t,\theta),
\end{equation} 
where the wave amplitude $\alpha$ remains undetermined to this order of approximation. It may be remarked here that the propagation along the other two families of characteristics can be handled in a similar manner. Using \eqref{alpha} in the $\rho^{(2)}_\theta,u^{(2)}_\theta$ and $p^{(2)}_\theta$ eliminant of \eqref{ordereps} we get the following transport equation for the wave amplitude:
\begin{equation} \label{main7}
\alpha_t + \alpha_x - \phi_t \left(\frac{\gamma +1}{2}\right)\alpha \alpha_\theta + \frac{j\alpha}{2x} =0.
\end{equation}
Furthermore, if we impose the condition that $\phi|_{x=1} = t,$ then $-\tfrac{\phi_t}{\phi_x} = 1$ yields 
\begin{equation}
\phi = t - (x-1).
\end{equation}
Accordingly, we may solve \eqref{main7} to obtain 
\begin{equation} \label{level}
\alpha x^{\tfrac{j}{2}} = \mbox{ constant },
\end{equation}
along the family of curves given by 
\begin{equation} \label{trajectories}
\frac{\,dx}{\,dt} = 1, \hspace{1cm} \frac{\,d\theta}{\,dt} = -\frac{\gamma +1}{2}\alpha.
\end{equation}
Let these curves (characteristic wavelets/or rays) be marked by the characteristic variable $\tau,$ which is chosen such that $t = \tau$ at the boundary $x=1;$ then integrating \eqref{trajectories} and using the definition of $\theta,$ we find that the equation of any \textit{wavelet,} $\tau = \mbox{ constant },$ is 
\begin{equation} \label{wavelet}
t - \tau = x-1 - \frac{\gamma +1}{2}ux^{\tfrac{j}{2}}J(x),
\end{equation}
where 
\begin{equation}
J(x) = \int_1^x r^{-\tfrac{j}{2}}\,dr = \left\{ \begin{array}{cc}
x-1, & \mbox{ plane }; \\
2(\sqrt{x} -1), & \mbox{ cylindrical }; \\
\log x, & \mbox{ spherical.}
\end{array}\right.
\end{equation}
In writing \eqref{wavelet}, we have replaced $\epsilon u^{(1)}$ by $u,$ as it is true in the leading order approximation. Also, equation \eqref{level} can be written as 
\begin{equation} \label{modsim2}
ux^{\tfrac{j}{2}} = v(\tau),
\end{equation}
where $v(\tau)$ is the value of $u$ on the boundary $x =1.$ Equations \eqref{wavelet} and \eqref{modsim2} constitute the desired modulated simple wave solution; infact, it follows that along a positive characteristic that intersects the boundary at time $\tau,$ we have $u = v(\tau),$ and the equation of the characteristic is 
\begin{equation} \label{char1}
t - \tau = x-1 - \frac{\gamma +1 }{2}v(\tau)J(x).
\end{equation}
Here, we envisage a weak shock wave that propagates into an undisturbed region with $u=0$ ahead of the shock. Thus, it follows that across the shock 
\begin{equation} \label{chareqns}
\dsp [u] = x^{-\tfrac{j}{2}}v(\tau_-),
\end{equation}
where $\tau_-$ is the wavelet immediately behind the shock. Also, if the shock path is $t=s(x),$ then the shock conditions \eqref{RH} implies that to the leading order approximation 
\begin{equation} \label{shockvelocity}
\frac{\,dt}{\,dx} = s'(x) = U^{-1} \sim 1 - \frac{\gamma + 1}{4}v(\tau_-)x^{-\tfrac{j}{2}},
\end{equation}
and the characteristic equation \eqref{wavelet} gives
\begin{equation} \label{shockeqn}
s(x) - \tau_- = x-1 - \frac{\gamma +1}{2} v(\tau_-)J(x).
\end{equation}
Equations \eqref{shockvelocity} and \eqref{shockeqn} determine the shock path $t = s(x)$ by elimintating the parameter $\tau_-.$ Further relations for the shock speed $U$ follow from \eqref{shockeqn}; this together with \eqref{shockvelocity} implies that on the shock
\begin{equation}
\frac{\gamma +1}{4}v(\tau_-)x^{-\tfrac{j}{2}} + \frac{\gamma +1}{2}\dot{v}(\tau)J(x) \frac{\,d\tau_-}{\,dx} = \frac{\,d\tau_-}{\,dx},
\end{equation}
which yields on integration 
\begin{equation} \label{49}
\frac{\gamma +1}{4} v^2(\tau_-)J(x) = \int_0^{\tau_-} v(s)\,ds.
\end{equation}
Also, since the characteristics $\tau_+ = 0$ ahead of the shock and the characteristic $\tau_-$ behind the shock meet at the shock at the same point and at the same time and since $v(\tau_+) = 0,$ it follows from \eqref{char1} that
\begin{equation} \label{derivative}
\lim_{\tau_- \to 0} \left( \frac{\tau_-}{v(\tau_-)}\right) = \frac{1}{\dot{v}(0)} = \frac{\gamma +1}{2}J(x).
\end{equation}
It may be noticed that \eqref{derivative} also follows from \eqref{49} in the limit $\tau_- \to 0,$ showing thereby that the integral in \eqref{49} is bounded at large distances (i.e., $x \to \infty$.) Thus, if the initial data is such that $v(\tau) \to 0$ as $\tau \to \tau_0,$ then some of the characteristics in the neighborhood of the zero wavelet do not overtake the shock from behind, and the asymptotic behaviour as $x \to \infty$ is 
\begin{equation} \label{51}
v(\tau_-) \sim \left(\frac{4b}{(\gamma +1)J(x)}\right)^{1/2},
\end{equation}
where $\left| \int_0^{\tau_0} v(\tau) \,d\tau\right| \leq b.$ Relation \eqref{51} together with \eqref{49} implies that 
\begin{equation}
[u] = u_- = v(\tau_-) x^{-\tfrac{j}{2}} \sim \left(\frac{4b}{(\gamma +1)J(x)}\right)^{1/2}x^{-\tfrac{j}{2}}.
\end{equation}
In other words, the shock decays like 
\begin{equation}
[u] \sim \left\{ \begin{array}{cc}
x^{-1/2}, & \mbox{ planar};\\
x^{-3/4}, & \mbox{ cylindrical };\\
x^{-1}(\log x)^{-\tfrac{1}{2}}, & \mbox{ spherical },
\end{array}\right.
\end{equation}
as $x \to \infty;$ this result is in agreement with that obtained in the preceding section. Furthermore, it follows from \eqref{char1} that between the shock and the limiting characteristic $\tau = \tau_0,$ we have
\begin{equation} \label{53}
x \sim t - \tau_0 + \frac{\gamma +1}{2} v(\tau)J(x),
\end{equation}
so that
\begin{equation} \label{54}
u = v(\tau)x^{-\tfrac{j}{2}} \sim \frac{2}{\gamma +1}(x - (t - \tau_0)) \left\{ \begin{array}{cc}
x^{-1}, & \mbox{ plane ;}\\
\frac{1}{2}x^{-1}, & \mbox{ cylinder; }\\
(x\log x)^{-1}, & \mbox{ spherical. }
\end{array}\right.
\end{equation}
The asymptotic behaviour of the first order discontinuity at the shock then follows from \eqref{54} by differentiating it partially with respect to $x,$ and then evaluating at the shock, namely
\begin{equation} \label{55}
u_x|_{-} = \frac{2}{\gamma +1} \{ K(x) + (x - s(x) + \tau_0) K'(x)\}; \hspace{1cm} K(x) =  \left\{ \begin{array}{cc}
x^{-1}, & \mbox{ plane ;}\\
\frac{1}{2}x^{-1}, & \mbox{ cylinder; }\\
(x\log x)^{-1}, & \mbox{ spherical. }
\end{array}\right.
\end{equation}
Evaluating \eqref{53} on the shock and then using it in \eqref{55} above, we find that to the leading order approximation
\begin{equation}
\dsp [u_x] \sim \frac{2}{\gamma + 1} \left\{ \begin{array}{cc}
x^{-1}, & \mbox{ plane; }\\
\frac{1}{2}x^{-1}, & \mbox{ cylindrical; }\\
(x\log x)^{-1}, & \mbox{ spherical, }
\end{array}\right.
\end{equation}
as $x \to \infty,$ which agrees with the decay behaviour outlined in the preceding section; indeed, it can be easily verified that in the weak shock limit $[u_x] \sim [p_x]$ to the leading order. 
\section{Simple Wave Formalism using Riemann Invariants}
System \eqref{1} can be written using vector matrix notation as 
\begin{equation} \label{55}
\mathbf{u}_t + A\mathbf{u}_x + B = 0,
\end{equation}
where $\mathbf{u} = \big(\begin{array}{ccc} \rho & u & p\end{array}\big)^{tr}, B = \big( \begin{array}{ccc} \tfrac{j\rho u}{x} & 0 & \tfrac{j\rho a^2 u}{x}\end{array}\big)^{tr}$ and $A$ is a $3 \times 3$ matrix that can be read off by an inspection of \eqref{1}. Since the eigenvalues of $A$ are $u \pm a$ and $u,$ and the corresponding left eigenvectors are $\big( \begin{array}{ccc} 0 & \pm a\rho & 1\end{array}\big)$ and $\big(\begin{array}{ccc}a^2 & 0 & -1\end{array}\big),$ respectively, it follows that equations \eqref{55} can be written in the characteristic form as
\begin{eqnarray}
\label{56} \frac{\,dp}{\,dt} \pm \rho a \frac{\,du}{\,dt} + \frac{j\rho u a^2}{x} = 0 \hspace{0.6cm} \mbox{ along } C_{\pm}:  \frac{\,dx}{\,dt} = u \pm a,\\
\label{57} \frac{\,dp}{\,dt} - a^2 \frac{\,d\rho}{\,dt} =0 \hspace{0.6cm} \mbox{ along } P: \frac{\,dx}{\,dt} = u,
\end{eqnarray}
where $C_{\pm}$ represent characteristics moving with velocity $\pm a$ relative to the velocity $u,$ and $P$ the particle path. Since $p = p(\rho, S)$ is a given function, where $S$ denotes the entropy, equation \eqref{57} shows that the entropy remains constant along each particle path $\tfrac{\,dx}{\,dt}= u.$ Further, as the jump in entropy at a shock is of the third order of smallness relative to the discontinuity in pressure, equation \eqref{57} implies that for a weak or even moderate strength shock, $S= S_0$ (constant) is a good approximation everywhere, and consequently equations \eqref{56} can be written as 
\begin{equation} \label{58}
\frac{\,d r^{\pm}}{\,dt} + \frac{j}{x}au = 0 \hspace{0.6cm} \mbox{ along } C^{\pm}: \frac{\,dx}{\,dt} = u \pm a,
\end{equation}
where $r^{\pm} = \left( \frac{2a}{\gamma -1} \pm u\right)$ are the Riemann invariants, and the differential operators on $r^{\pm}$ are just the operators of differentiation along the characteristics $C^{\pm}.$ The flow described by the solution of \eqref{58} is a simple wave; in the plane $(j=0)$ case, equation (57b) is dispensed with quickly by arguing that $r^-$ is constant everywhere. For a nonplanar case, the situation is different. As the jump in $r^-$ across a shock is also of the third order of smallness relative to the jump $[p],$ and the integral of the last term in (57b) taken along $C^-$ is small keeping in view that the simple wave argument applies in general near the wave front of any disturbance propagating into a uniform state at rest, it follows that the constancy of entropy $S$ and the Riemann invariant $r^-,$ i.e., $S= S_0$ and $r^- = \frac{2}{\gamma-1}$ are a good approximation to the first order, and thus the first equation in \eqref{56} provides a single first order equation for $u,$ namely,
\begin{equation} \label{59}
u_t + \left( 1 + \frac{\gamma +1}{2}u \right) u_x + \frac{j}{x}u\left( 1 + \frac{\gamma -1}{2}u\right) = 0.
\end{equation}
Equation \eqref{59} can be solved exactly to yield 
\begin{equation} \label{60}
u\left( 1 + \frac{\gamma -1}{2}u\right)^{\frac{2}{\gamma-1}} = v(\tau)x^{-\tfrac{j}{2}},
\end{equation}
where $\tau$ is the characteristic variable to be determined from $\tfrac{\partial t}{\partial x} = \big(1 + \tfrac{\gamma+1}{2}u\big)^{-1} \sim 1 - \frac{\gamma+1}{2}u.$

Thus, a uniformly valid approximation at the leading order is
\begin{equation}
\begin{array}{l}
u = v(\tau)x^{-\tfrac{j}{2}}, \\
t = \tau + (x-1) - \frac{\gamma +1}{2}v(\tau)J(x), 
\end{array}
\end{equation}
which are exactly the same as (39)-(41), and hence the evolutionary behaviour of shock and the first order discontinuity turn out to be the same as discussed in the preceding sections.
\section{Method of Relatively Undistorted Waves}
This method proposed by Varley and Cumberbatch \cite{varley, varley1} is based on a scheme of successive approximations to a system of hyperbolic equations. The theory does not make any assumptions on the amplitude of the disturbance and yields exact results for simple waves, acceleration waves, and thus for the formation of shock fronts. Here, we show that this method yields the same results as those derived in section 7 for the evolution of shock strength and the accompanying first order discontinuity. 

Under the transformation of coordinates $(x,t) \mapsto (x,\tau),$ where $\tau(x,t)$ is a wavelet defined by $t= t(x,\tau),$ and any vector valued function $\mathbf{u}(x,t)$ transforms as $\mathbf{u}(x,t) \mapsto \mathbf{U}(x,\tau),$ system \eqref{55} may be written as 
\begin{equation} \label{62}
(A - t_x^{-1}\mathcal{I})\mathbf{u}_t = t_x^{-1}(A\mathbf{U}_x + B),
\end{equation}
where $\mathcal{I}$ is the identity matrix of order $3.$ A solution vector $\mathbf{u}$ is said to be relatively undistorted if 
\begin{equation} \label{63}
\|\mathbf{U}_x\| << \|\mathbf{u}_x\|,
\end{equation}
implying thereby that $\mathbf{u}_x \sim t_x \mathbf{u}_t.$ Equations \eqref{62} and \eqref{63} are compatible if $\|B\| = \OO(1)\|A\mathbf{U}_x\|,$ while $t_x^{-1}$ is an eigenvalue of $A;$ this implies that the wavelets are the characteristic curves of \eqref{55}. Subsequently it follows from \eqref{62} that $\mathbf{U}$ must satisfy the compatibility condition 
\begin{equation}\label{64}
\mathbf{L}(A\mathbf{U}_x + B) = 0,
\end{equation}
where $\mathbf{L} = \big( \begin{array}{ccc} 0 & \rho a & 1\end{array}\big)$ is the left eigenvector of $A$ corresponding to the eigenvalue $t_x^{-1} = u+a.$ We look for the solution of \eqref{55} in the region $ x >1,$ where the motion associated with the eigen-mode $t_x^{-1}= u+a$ is excited by the boundary condition 
\begin{equation}  \label{65}u(1,\tau) = v(\tau),
\end{equation}
where $v$ is a smooth bounded function, i.e., $|v| = \OO(1).$ We let $\tau = $ constant denote the characteristic curve defined by $t_x^{-1} = u+a,$ and parametrize $\tau$ by defining $\tau = t$ at $x=1.$ Thus, equation \eqref{62} to a first approximation yields 
\begin{equation} \label{66}
(A- t_x^{-1} \mathcal{I}) \mathbf{U}_\tau = 0,
\end{equation}
which together with \eqref{64} and $t_x^{-1} = u+a$ enables us to determine the unknowns $\mathbf{U}$ and $t(x,\tau).$ Equation \eqref{66} implies that $\mathbf{U}$ is collinear to the right eigenvector $\mathbf{R} = \big(\begin{array}{ccc} \rho & a & \rho a^2\end{array}\big)$ of $A$ corresponding to the eigenvalue $u+a.$ For a wave moving into an undisturbed region, equation \eqref{66} integrates to give 
\begin{equation} \label{67}
\rho = \left( 1 + \frac{\gamma -1}{2}u\right)^{\tfrac{2}{\gamma-1}}, \hspace{0.6cm} p = \frac{1}{\gamma}\left(1 + \frac{\gamma-1}{2}u\right)^{\tfrac{2\gamma}{\gamma-1}}, \hspace{0.6cm} a = 1 + \frac{\gamma-1}{2}u. 
\end{equation}
Equations \eqref{67} hold throughout the region at any $x$ on any wavelet $\tau = $ constant; these may be used in \eqref{64} to yield the following transport equation for $u,$ i.e.,
\begin{equation} \nonumber
\left(1 + \frac{\gamma+1}{2}u \right)u_x + \frac{j}{2x}u\left(1 + \frac{\gamma-1}{2}u \right) = 0,
\end{equation}
which on integration yields
\begin{equation} \nonumber
u \left( 1 + \frac{\gamma-1}{2}u\right)^{\tfrac{2}{\gamma-1}} = v(\tau)x^{\tfrac{-j}{2}},
\end{equation}
where $\tau$ is the characteristic variable to be determined from $t_x^{-1} = \big( 1+ \tfrac{\gamma+1}{2}u\big)^{-1} \sim \big( 1- \tfrac{\gamma+1}{2}u\big);$ this indeed, is precisely the solution (61). Thus, the relatively undistorted wave approximation is a consequence of the simple wave formalism using Riemann invariants. 

At this juncture, it is worthwhile to remark that the Chester-Chisnell-Whitham (CCW) approximation, which consists in applying the shock conditions \eqref{RH} to the exact characteristic equation \eqref{56} along $C_+$, yields the following transport equation for the shock
\begin{equation} \label{68}
\frac{Ug(U)}{U^2-1} \frac{\,dU}{\,dx} + \frac{j}{x} = 0,
\end{equation}
where $\tfrac{\,d}{\,dx}$ denotes the derivative along $C_+: \tfrac{\,dx}{\,dt} = (u+a)|_{\mbox{\small{shock}}},$ and $g(U) = \big(1 + 2\big(\tfrac{\mu}{\nu}\big)^{1/2} + U^{-2}\big)\{1 + \tfrac{U^2-1}{(\mu\nu)^{1/2}}\}.$ It is indeed interesting to note that the transport equation (8) for the shock, in view of relation \eqref{RH}, assumes the following form by neglecting the precursor disturbances from the rearward flow (i.e., on setting $[p_x]=0$)
\begin{equation} \label{69}
\frac{UG(U)}{U^2-1} \frac{\,dU}{\,dx} + \frac{j}{x}=0,
\end{equation}
where 
$G(U) = (\gamma +1) \{2U^2 \nu^{-1} + \mu^{-1} (U^2+1)\};$ it may be noticed that \eqref{68} and \eqref{69}, except for the functional forms of $g(U)$ and $G(U)$ and the differential operators, share a close structural resemblance. Indeed, it is remarkable that for weak shocks, $U \sim 1,$ both $G(U)$ and $g(U) \sim 4,$ and the evolutionary behaviour governed by \eqref{68} and \eqref{69} coincide. 

\section{Results and Conclusion}
Except some impirical methods, which have been developed in the past, no analytical
method exists for describing the evolutionary behaviour of a shock wave without
limiting its strength. In the present paper, we generate a completely intrinsic description
of plane, cylindrical and spherical shock waves of weak or even moderate strength,
propagating into a uniform flow of a polytropic gas at rest. We use four different approximation methods to study the asymptotic evolutionary behaviour of shocks and
the associated first order jump discontinuity, which may be interpreted as the disturbance that catches up with the shock from behind. Singular surface theory yields a set
of coupled equations for the shock strength and the higher order jump discontinuities
which couple the shock motion with the rearward flow. It is noticed that the transport
equations for the shock strength $[p]$ and the first order jump discontinuity $[p_x]$ are nonlinear in $[p]$ and $[p_x]$ respectively; however, the transport equations for $[p_{xx}]$ is linear
in $[p_{xx}]$, and this is true for all other higher order jump discontinuities. The transport
equation for the shock strength shows that the evolutionary behaviour of the shock
at any instant is influenced by the shock curvature and the first order jump discontinuity,
representing the effect of precursor disturbances that overtake the shock from
behind. This transport equation, in some sense, generalizes the CCW approximation
\cite{whitham}, which although does not account for the effects of the disturbances that overtake
the shock from behind, yet its success is unexpectedly remarkable for a particular
class of problems. We show in the preceding section that the transport equation for the
shock strength, in the limit of a weak shock and vanishing $[p_x]$, reduces exactly to the
one obtained by the CCW approximation. The truncated system of transport equations for $[p]$ and $[p_x]$ is solved for uniformly valid solutions using a singular perturbation technique \cite{nayfeh} that eliminates secular terms in the perturbed solution; the shock strength is assumed small, but the accompanying first order jump discontinuity may be finite (i.e., $\OO(1)$) or small. It is shown that the precursor disturbances evolve to the leading order like an acceleration wave. Indeed, the decay of the induced discontinuity is much faster than the decay of the shock; it is noticed that for the case $[p_x] = \OO(\epsilon),$ the shock decays faster as compared to the case when $[p_x] = \OO(1).$ Asymptotic results for the decay of weak shocks fully agree with those given in \cite{whitham,vdsbook}. It is shown that the approximation methods of weakly nonlinear geometrical optics, simple waves involving  Riemann invariants \cite{sirovich}, and the relatively undistorted wave approximation yield the same results within their common range of validity. 

In conclusion, we observe that an analytical description of  the evolutionary behaviour of shocks of arbitrary strength in an ideal gas remains an open problem because the truncation approximation, used in rendering an infinite hierarchy of transport equations as a closed finite system, may not yield a uniformly valid solution for all running lengths. The present work envisages four different methods to study the evolutionary behaviour of weak shocks and the accompanying first order discontinuity; a thorough comparative analysis of these methods has been presented, which to the best our knowledge has not been done earlier. Indeed, the singular surface theory gives much more powerful results then does the method of characteristics that forms the basis of rest of the approximation methods used herein. The truncation approximation leading to the closed system consisting of  a coupled pair of transport equations for $[p]$ and $[p_x]$ can be regarded as a good approximation for the hierarchy of the system governing shock evolution. Further, as the relatively undistorted wave approximation, and hence the method of weakly nonlinear geometrical optics, which satisfies the condition needed for an approximation to be relatively undistorted, lead to the first order the same evolution laws for the shock and the associated first order discontinuity, it transpires that both these approximation methods are a consequence of the simple wave formalism using Riemann invariants. However, the relatively undistorted approximation makes no assumption on the amplitude of the disturbance, unlike the later which is a manifestation of the method of multiple scales. Lastly, a comparison with the CCW approximation reveals that the singular surface theory  yields the shock-transport equation which can be regarded as a generalization of the CCW approximation.

\section*{Acknowledgements}
The authors are thankful to the referee for his elaborate suggestions in improving the quality of the manuscript, and for pointing out the reference \cite{na}. RV is grateful to the KVPY Cell, Department of Science and Technology, India, for an undergraduate fellowship, and to Indian Institutes of Technology Bombay and Gandhinagar for their hospitality during his visit.
\bibliographystyle{amsplain}

\end{document}